\documentclass[12pt,a4paper,reqno]{amsart}

\usepackage[latin1]{inputenc}     
\usepackage[T1]{fontenc}          
\usepackage[english]{babel}        

\usepackage{amsmath,amsfonts,amssymb,amsthm} 
\usepackage{mathtools}            
\usepackage{mathrsfs}
\usepackage{physics}

\usepackage{url}
\usepackage{cite}
\usepackage{enumitem}

\usepackage{microtype}            
\usepackage{enumitem}             
\usepackage{graphicx}             
\usepackage{booktabs}             
\usepackage{geometry}             
\usepackage[colorlinks=true, linkcolor=blue, citecolor=blue, urlcolor=blue, pdfstartview=FitH]{hyperref}

\newtheorem{theorem}{Theorem}

\newtheorem{lemma}{Lemma}

\title[Stabilization of an EBB with damping and interior time delay]{Stabilization of an Euler--Bernoulli beam equation with distributed damping and interior time delay}

\author[B. B. J. Siriki]{Ben B. J. Siriki}
\address{Universit\'e Nangui Abrogoua, Abidjan, C\^ote d'Ivoire}
\email{\href{mailto: benbjsiriki@gmail.com}{benbjsiriki@gmail.com}}

\author[N. F. Alla]{N'Guessan F. Alla}
\address{Universit\'e F\'elix  Houphou\"et--Boigny,
	Abidjan, C\^ote d'Ivoire.}
\email{allafabrice.nguessan@ufhb.edu.ci}

\author[A. Coulibaly]{Adama Coulibaly}
\address{Universit\'e F\'elix  Houphou\"et--Boigny,
	Abidjan, C\^ote d'Ivoire}
\email{couliba@yahoo.fr}

\begin{document}
	
	\begin{abstract}
		In this paper, we investigate the effect of a fixed time delay on the internal stabilization of an Euler--Bernoulli beam clamped at one end with distributed damping. Under appropriate assumptions on the feedback weights, we first establish that the problem is well-posed using semigroup theory. Furthermore, using the multiplier method combined with a suitable Lyapunov functional, we prove the exponential stability of the system.
	\end{abstract}
	
	\maketitle
	
	\tableofcontents

\section*{Introduction}

This study aims to analyze the impact of a fixed internal time delay on the internal stabilization of a uniform Euler--Bernoulli beam, characterized by its mass per unit length $m$ and flexural rigidity $EI>0$. The vibratory dynamics of the system are governed by the following initial-boundary value problem:
\begin{equation} \label{pde-system}
	\left\{ 
	\begin{aligned}
		& m y_{tt} + EI y_{xxxx} + u_\tau(x,t) = 0, \quad x \in (0,\ell), \, t>0, \\
		& y(0,t) = y_x(0,t) = y_{xx}(\ell,t) = y_{xxx}(\ell,t) = 0, \quad t>0,  \\
		& y(x,0) = y_0(x), \, y_t(x,0) = y_1(x), \quad x \in [0, \ell], \\
		& y_t(x,t-\tau) = f_0(x,t-\tau), \quad x \in (0, \ell), \, t \in (0, \tau),
	\end{aligned} 
	\right. 
\end{equation}
where $y(x,t)$ represents the deflection of the beam at position $x$ and time $t$, while $u_\tau(x,t)$ denotes the internal control law, which depends on a fixed time delay $\tau>0$. The boundary condition $\eqref{pde-system}_2$ means that the beam is clamped at its end $x=0$ and free at the opposite end $x=\ell$. Finally, $(y_0, y_1)$ is the initial state of the beam, and $f_0$ indicates the memory of controller.

The stabilization of Euler--Bernoulli beam equation is central to control theory, driven by both its engineering applications and its mathematical significance \cite{Du2022, Chakravarthy2010, Balakrishnan1987, He2014, Wright2003}. This challenge intensifies significantly when time delays are introduced, as even a small latency in the feedback loop can trigger instability \cite{Gumowski1968, Datko1985, Datko1991, Datko1993}. For this reason, designing robust control strategies to neutralize delay-induced instabilities remains a top priority in this field.

In this perspective, two main control strategies are widely distinguished in the literature for time-delay Euler--Bernoulli beam equations, depending on where the control action is applied.
On the one hand, boundary control consists in applying feedback directly at the boundaries of the structure. In their pioneering work \cite{Guo2009}, Guo et al. considered an output feedback stabilization problem for the system:
\begin{equation} \label{pde-r1}
	\left\{ 
	\begin{aligned}
		& w_{tt} + w_{xxxx} = 0, \quad 0 < x < 1, \, t>0,\\
		& w(0,t) = w_x(0,t) = w_{xx}(1,t) = 0,  \\
		& w_{xxx}(1,t) = u(t), \\
		& y(t) = w_t(1,t-\tau), \quad t \geqslant \tau, \\
		& w(x,0) = w_0(x), \, w_t(x,0) = w_1(x), \\
	\end{aligned} 
	\right. 
\end{equation}
where $u$ is the control, $(w_0,w_1)$ the initial state, and $y$ the delayed measurement. The authors exponentially stabilize the system using the control law
\begin{gather*}
	u(t) = 
	\left\{ 
	\begin{aligned}
		& k_2 \widetilde{w}_t(1,t) = k_2 \widehat{w}_s(1,t,t), \quad t > \tau,
		\\
		& 0, \quad t \in [0,\tau],
	\end{aligned} 
	\right.
\end{gather*}
designed via an observer-predictor scheme based on the estimated state $(\tilde{w}, \tilde{w}_t)$ with $k_2>0$. Subsequently, numerous contributions have deepened this direction, proposing various boundary control laws to effectively suppress the vibrations of such structures. A detailed review of these approaches and their recent developments can be found in \cite{Cao2026, Mahfoudhi2026, Azzaoui2026, Salhi2026, Camasta2025, Rahmtati2022, Siriki2025, Siriki2026, Siriki2026a, Lhachemi2019}. \newline
On the other hand, the control action can be distributed inside the spatial domain: this is referred to as internal (or distributed) control. The internal stabilization of Euler--Bernoulli beam models was pioneered by Feng et al. In \cite{Feng2019}, the state-feedback stabilization of the following system was addressed:
\begin{equation} \label{pde-r2}
	\left\{ 
	\begin{aligned}
		& w_{tt}(x,t) + w_{xxxx}(x,t) + u(x,t-\tau) = 0, \quad x \in (0,1), \, t>0, \\
		& w(0,t) = w_{x}(1,t) = 0,  \\
		& w_{xx}(0,t) = w_{xxx}(1,t) = 0, \\
		& w(x,0) = w_0(x), \, w_t(x,0) = w_1(x), \quad x \in (0,1)\\
		& u(x,s-\tau) = h(x,s), \quad x \in (0,1), \, s \in (0,\tau),		
	\end{aligned} 
	\right. 
\end{equation}
where $h(x,s)$ stands for the memory function of the controller. Employing the backstepping method \cite{Smyshlyaev2009}, they designed a feedback control law
\begin{equation}
	\begin{aligned}
		u(x,t) 
		& = 
		\int_0^\tau \int_0^1 q(x,\tau;y,r) z(y,r,t)\, dy dr 
		+ \int_0^1 \gamma(x,\tau,r) w(y,t)\, dy 
		\\
		& \quad 
		+ \int_0^1 \eta(x,\tau,y) w_t(y,t)\, dy, \quad x \in (0,1), \, t \geqslant 0,
	\end{aligned} 
\end{equation} 
via a bounded, linear, and invertible integral operator-structured using appropriate kernel functions such as $q(x,\tau;y,r)$, $\gamma(x,\tau,r)$, and $\eta(x,\tau,y)$, which maps the original plant into an exponentially stable target system. 

However, internal feedback mechanisms remain considerably less explored than their boundary counterparts. Yet, distributed controllers provide greater flexibility and allow for a more effective mitigation of localized physical phenomena \cite{Fanou2026}, such as periodic disturbances, structural resonances, or spatial constraints. 
To achieve the exponential stabilization of system \eqref{pde-system}, we propose the following localized internal feedback control law: 
\begin{equation} \label{control-law}
	u_\tau(x,t) = a(x) \left( \mu_1 y_t(x,t) + \mu_2 y_t(x,t-\tau) \right), \quad x \in (0,\ell), \, t \geqslant 0.
\end{equation}
inspired by the framework of Nicaise and Pignotti \cite{Nicaise2006}, where $\mu_1 > 0$, $\mu_2 \neq 0$, and $a \in L^\infty(0,\ell)$ is a function satisfying
\begin{equation} \label{damper-cdt}
	\left\{ 
	\begin{aligned}
		& a(x) \geqslant 0 \quad \text{ a.e. in } (0,\ell),
		\\
		& a(x) > a_0 > 0 \quad \text{ a.e. in } \omega,
	\end{aligned} 
	\right.
\end{equation}
with $\omega \subset (0,\ell)$ being a non-empty open neighbourhood of $\ell$.
To the best of our knowledge, the internal stabilization of an Euler--Bernoulli beam under such a localized, delayed distributed feedback law remains entirely open in the literature. First, we reformulate system \eqref{pde-system} by introducing a new state variable to handle the time delay. Then, under the dominance condition
\begin{gather} \label{internal-coef-cdt}
	|\mu_2| < \mu_1,
\end{gather} 
we first prove that the problem is well-posed by invoking the semigroup theory \cite{Pazy1983}. Secondly, we establish the exponential decay of the system's energy via a suitable Lyapunov functional coupled with the multiplier technique, departing from the observability-based approach of \cite{Nicaise2006}.

The rest of this paper is organized as follows. Section \ref{sec1} is devoted to preliminary reminders and tools necessary for establishing our main results. In Section \ref{sec2}, we study the well-posedness of the problem \eqref{pde-system} using semigroup theory. Finally, Section \ref{sec3} is entirely dedicated to proving the exponential stability of the system.

\section{Preliminaries} \label{sec1}
This section gathers the essential tools for our main analysis.
Throughout this paper, $C$ denotes a generic constant.
%

\subsection{Augmented system}

To handle the time delay, we follow \cite{Nicaise2006} and introduce the auxiliary variable $z(x, \rho, t)$ defined as
\begin{equation}
	z(x,s,t) := y_t(x,t-s\tau), \quad x \in (0,\ell), \, s\in (0,1), \, t>0.
\end{equation} 
By direct computation, $z$ satisfies the transport equation:
\begin{equation}
	\left\{ 
	\begin{aligned}
		& \tau z_t + z_s = 0, \quad x \in (0,\ell), \, s\in (0,1), \, t>0, \\
		& z(x,0,t) = y_t(x,t), \quad x \in (0,\ell), \, t>0, \\
		& z(x,1,t) = y_t(x,t-\tau), \quad x \in (0,\ell), \, s\in (0,1), \, t>0, \\
		& z(x,s,0) = f_0(x,-s\tau), \quad x \in (0, \ell), \, s \in (0, 1).
	\end{aligned} 
	\right.
\end{equation}
Consequently, system \eqref{pde-system} takes the augmented form:
\begin{equation} \label{pde-aug-system}
	\left\{ 
	\begin{aligned}
		& m y_{tt} + EI y_{xxxx} + a(x) \big( \mu_1 y_t(x,t) + \mu_2 \, z(x,1,t) \big) = 0, \\
		&\hspace{6.5cm} x \in (0,\ell), \, t>0, \\
		& \tau z_t + z_s = 0, \quad x \in (0,\ell), \, s\in (0,1), \, t>0, \\
		& y(0,t) = y_x(0,t) = y_{xx}(\ell,t) = y_{xxx}(\ell,t) = 0, \quad t>0,  \\
		& z(x,0,t) = y_t(x,t), \quad x \in (0,\ell), \, t>0, \\
		& y(x,0) = y_0(x), \, y_t(x,0) = y_1(x), \quad x \in [0, \ell], \\
		& z(x,s,0) = f_0(x,-s\tau), \quad x \in (0, \ell), \, s \in (0, 1).
	\end{aligned} 
	\right.
\end{equation} 

\subsection{Energy of the system}

Let $E(t)$ be the energy of system \eqref{pde-system} defined by:
\begin{gather} \label{energy-system}
	\begin{aligned}
		E(t) & := \frac{m}{2} \int_{0}^{\ell} \, y_{t}^2 \, dx + \frac{EI}{2} \int_{0}^{\ell} \, y_{xx}^2 \, dx 
		+ \frac{\xi}{2} \int_{0}^{\ell} a(x) \int_{0}^{1} \, y_t^2(x,t-s\tau) \, ds \, dx,
	\end{aligned}
\end{gather}
where $\xi$ is a positive constant. 

\begin{lemma} \label{l1a}
	Assume that \eqref{internal-coef-cdt} holds. If the parameter $\xi$ satisfies the condition
	\begin{gather} \label{aux-coef-cdt}
		\tau |\mu_2| \leqslant \xi \leqslant \tau \big( 2 \mu_1 - |\mu_2| \big)
	\end{gather}
	then, the energy $E(t)$ is non-increasing over time, i.e.:
	\begin{align} \label{decrease-energy}
		\begin{split}
			\frac{d}{dt} E(t) 
			& \leqslant 
			-\left( \mu_1 - \frac{|\mu_2|}{2} - \frac{\xi \tau^{-1}}{2} \right) \int_{0}^{\ell} a(x) y_t^2 \, dx 
			\\
			& \quad - \left( -\frac{|\mu_2|}{2} + \frac{\xi \tau^{-1} }{2}\right) \int_{0}^{\ell} a(x) y_t^2(x,t-\tau) \, dx.
		\end{split}
	\end{align}
\end{lemma}

\begin{proof}
	First, multiplying \eqref{pde-aug-system}$_1$ by $y_t$ and integrating over $(0,\ell)$ with respect to $x$, we obtain
	\begin{equation} \label{proof-l1a-1}
		\begin{split}
			\int_{0}^{\ell} \, \Big( m y_{tt} + EI y_{xxxx} + a(x) \big( \mu_1 y_t(x,t) + \mu_2 z(x,1,t) \big) \, dx \Big) y_t \, dx = 0.
		\end{split}
	\end{equation}
	Integrating by parts, we get
	\begin{gather} 
		\begin{aligned}[b]
			\int_{0}^{\ell} \, y_{tt} y_t \, dx
			&= 
			\frac{d}{dt} \left(
			\frac{1}{2} \int_{0}^{\ell} \, y_{t}^2 \, dx \right)
		\end{aligned} 
		\\
		\begin{aligned}[b]
			\int_{0}^{\ell} \, y_{xxxx} y_t \, dx 
			&= 
			\big[y_{xxx} y_t - y_{xx} y_{tx} \big]_{x=0}^{x=\ell} + \int_{0}^{\ell} \, y_{xx} y_{txx} \, dx 
			\\
			&=
			\frac{d}{dt} \left(
			\frac{1}{2} \int_{0}^{\ell} \, y_{xx}^2 \, dx \right),
		\end{aligned}
	\end{gather}
	where we used the boundary condition \eqref{pde-aug-system}$_3$. Therefore, relation \eqref{proof-l1a-1} simplifies to
	\begin{equation} \label{proof-l1a-2}
		\begin{split}
			& \frac{d}{dt} \left(
			\frac{m}{2} \int_{0}^{\ell} \, y_{t}^2 \, dx + \frac{EI}{2} \int_{0}^{\ell} \, y_{xx}^2 \, dx \right)
			\\
			& \quad  = - \mu_1 \int_{0}^{\ell} a(x) y_t^2 \, dx 
			- \mu_2 \int_{0}^{\ell} a(x) z(x,1,t) y_t \, dx.
		\end{split}
	\end{equation}
	Next, equation $\eqref{pde-aug-system}_2$ implies
	\begin{gather} \label{proof-l1a-3}
		\xi \int_{0}^{\ell} a(x) \int_{0}^{1} \, z_t z \, ds \, dx= - \xi \tau^{-1} \int_{0}^{\ell} a(x) \int_{0}^{1} \, z_s z \, ds \, dx.
	\end{gather}
	An integration by parts with respect to $s$, together with the boundary condition $\eqref{pde-aug-system}_4$, yields
	\begin{align}
		\int_{0}^{1} z_t z \, ds 
		& =  \frac{d}{dt} \left( \frac{1}{2} \int_{0}^1 z^2  \, ds \right); \label{proof-l1a-4}
		\\
		\int_{0}^{1} \, z_s z \, ds 
		& = 
		\frac{1}{2} \Big( z(x,1,t)^2 - y_t^2 \Big). \label{proof-l1a-6}
	\end{align}
	Combining identities \eqref{proof-l1a-4}--\eqref{proof-l1a-6}, it follows that
	\begin{gather} \label{proof-l1a-7}
		\begin{aligned}[b]
			\frac{d}{dt} \left( \xi \int_{0}^{\ell} a(x) \int_{0}^{1} \, z^2 \, ds \, dx \right) 
			&=
			\xi \tau \int_{0}^{\ell} a(x) \Big( y_t^2 - z(x,1,t)^2 \Big) \, dx.
		\end{aligned}
	\end{gather}
	Finally, adding \eqref{proof-l1a-2} to \eqref{proof-l1a-7} leads to
	\begin{equation} \label{proof-l1a-8}
		\begin{aligned}
			& \frac{d}{dt} \left(
			\frac{m}{2} \int_{0}^{\ell} \, y_{t}^2 \, dx + \frac{EI}{2}  \int_{0}^{\ell} \, y_{xx}^2 \, dx 
			+ \frac{\xi}{2} \int_{0}^{\ell} a(x) \int_{0}^{1} \, z^2 \, ds \, dx \right) 
			\\ 
			& = 
			- \mu_1 \int_{0}^{\ell} a(x) y_t^2 \, dx
			+ \frac{\xi \tau^{-1}}{2} \int_{0}^{\ell} a(x) \big( y_t^2 - z(x,1,t)^2 \big) \, dx 
			\\
			& \quad - \mu_2 \int_{0}^{\ell} a(x) z(x,1,t) y_t \, dx.
		\end{aligned}
	\end{equation}
	In addition, by Young's inequality, it follows that
	\begin{align} \label{proof-l1a-9}
		\big| \mu_2 z(x,1,t) y_t \big|
		& \leqslant \frac{|\mu_2|}{2} \Big( z(x,1,t)^2 + y_t^2 \Big).
	\end{align}
	Substituting estimate \eqref{proof-l1a-9} into \eqref{proof-l1a-8} directly implies \eqref{decrease-energy}. The non-positivity of the right-hand side is guaranteed whenever $\xi$ satisfies \eqref{aux-coef-cdt}.
\end{proof}

\section{Existence and uniqueness of solution} \label{sec2}

This section is devoted to the well-posedness of problem \eqref{pde-system} via semigroup theory. Specifically, we recast the augmented system \eqref{pde-aug-system} into an abstract Cauchy problem and apply the Lumer--Phillips theorem.

\subsection{Abstract formulation}

Let $H^2_e(0,\ell) := \{y \in H^2(0,\ell) \mid y(0) = y'(0) = 0\}$. We define the state space $\mathscr{H}$ as
\begin{align}
	\mathscr{H} := H^2_e(0,\ell) \times L^2(0, \ell) \times L^2\big((0,\ell) \times (0,1)\big),
\end{align}
equipped with the inner product
\begin{align}
	\langle U_1, U_2 \rangle 
	= 
	m \int_{0}^{\ell} v_1 \overline{v_2} \, dx + EI \int_{0}^{\ell} \, y_1'' \overline{y_2''} \, dx + \xi \int_{0}^{\ell} a(x) \int_{0}^{1} z_1 \overline{z_2}\, ds \, dx,
\end{align}
for all $U_i = (y_i, v_i, z_i) \in \mathscr{H}$, $i = 1, 2$. It is straightforward to check that $\mathscr{H}$ is a Hilbert space. 
Next, on the space $\mathscr{H}$, we define the linear operator $\mathcal{A}$ with domain $\mathscr{D}(\mathcal{A}) \subset \mathscr{H}$ given by
\begin{gather}
	\begin{gathered}
		\mathscr{D}\big(\mathcal{A}\big) := 
		\left\{
		U \in \mathscr{H}: \, 
		\begin{matrix}
			y \in H^4(0,\ell), \, v \in H^2_e(0,\ell), 
			\\ 
			z \in L^2\big( 0,\ell; H^1(0,1) \big),
			\\
			z(\cdot,0) = v, \, y''(\ell) = y^{(3)}(\ell) = 0
		\end{matrix}
		\right\}
	\end{gathered}
	\\
	\mathcal{A} U :=
	\begin{pmatrix}
		v \\ -\frac{EI}{m} y^{(4)} - \frac{a}{m} \Big( \mu_1 v + \mu_2 z(\cdot,1) \Big) \\ -\tau^{-1} z_s
	\end{pmatrix}.
\end{gather}
With these definitions at hand, the system \eqref{pde-aug-system} can be written into the abstract Cauchy problem
\begin{gather} \label{Cauchy-pb}
	\left\{
	\begin{aligned}
		& \frac{d}{dt} U(t) = \mathcal{A} U(t), 
		\\
		& U(0) = U_0 = \big( y_0, \, y_1, \, f_0 \big).
	\end{aligned}
	\right.
\end{gather}

\subsection{Well-posedness of system}

Before establishing the well-posedness of system \eqref{pde-system}, we prove the following lemma regarding the norm equivalence on $H^2_e(0,\ell)$.

\begin{lemma} \label{l1}
	The mapping $|\cdot|_2: \, H^2_e(0,\ell) \to \mathbb{R}$ defined by
	\begin{equation} \label{norm:def}
		|w|_2 := \left( \int_0^\ell \big| w''(x) \big|^2 \, dx +   \int_0^\ell |w(x)|^2 \, dx \right)^{\frac{1}{2}}, \; \forall w \in H^2_e(0,\ell).
	\end{equation}
	defines a norm on $H^2_e(0,\ell)$ equivalent to the standard $H^2(0,\ell)$--norm. That is, there exist two positive constants $C_1, \, C_2$ such that
	\begin{gather} \label{norm:equiv}
		C_1 \lVert w \rVert_{H^2(0,\ell)} \leqslant |w|_2 \leqslant C_2 \lVert w \rVert_{H^2(0,\ell)}, \quad \forall w \in H^2_e(0,\ell).
	\end{gather}
\end{lemma}

\begin{proof}
	The verification that $\lvert \cdot \rvert_2$ satisfies the norm axioms on $H^2_e(0,\ell)$ is standard and omitted here. 
	Next, we prove the equivalence property. Let $w \in H^2_e(0,\ell)$. By definition, we immediately have $\lvert w \rvert_2 \leqslant \lVert w \rVert_{H^2(0,\ell)}$. 
	Furthermore, using the triangle and Cauchy--Schwarz inequalities, we get
	\begin{gather} 
		\begin{aligned}[b]
			\int_0^\ell |w'(x)|^2 \, dx 
			& \leqslant 
			\int_0^\ell x \int_0^x  |w''(\sigma)|^2 \, d\sigma \, dx
			\\
			& \leqslant \frac{\ell^2}{2} \lVert w'' \rVert^2_{L^2(0,\ell)}.
		\end{aligned}
	\end{gather}
	Consequently, we arrive at
	\begin{gather}
		\lVert w \rVert_{H^2(0,\ell)}^2 
		\leqslant \left( 1 + \frac{\ell^2}{2} \right) |w|_2^2,
	\end{gather}
	which completes the proof.
\end{proof}

\begin{theorem}
	Assume that conditions \eqref{internal-coef-cdt} and \eqref{aux-coef-cdt} holds. For any $U_0 \in \mathscr{H}$, problem \eqref{pde-aug-system} admits a unique mild solution $U \in C\left(\mathbb{R}_+, \mathscr{H}\right)$.
	Moreover, if $U_0 \in \mathscr{D}(\mathcal{A})$, then $U$ is a classical solution satisfying $U \in C(\mathbb{R}_+, \mathscr{D}(\mathcal{A})) \cap C^1(\mathbb{R}_+, \mathscr{H})$.
\end{theorem}

\begin{proof}
	First, we prove that the operator $\mathcal{A}$ is dissipative. Let $U \in \mathscr{D}( \mathcal{A} )$. By definition of the inner product in $\mathscr{H}$, we have
	\begin{equation} \label{proof-t1-1}
		\begin{aligned}
			\langle \mathcal{A}U, U \rangle 
			&=
			- \mu_1 \int_{0}^{\ell} a(x) |v(x)|^2 \, dx - EI \int_{0}^{\ell} y^{(4)} \overline{v} \, dx + EI \int_{0}^{\ell} v'' \overline{y''} \, dx 
			\\
			& \quad - \mu_2 \int_{0}^{\ell} a(x) z(x,1) \overline{v} \, dx  - \xi \tau^{-1} \int_{0}^{\ell} a(x) \int_{0}^{1} z_s \overline{z}\, ds \, dx.
		\end{aligned}
	\end{equation}
	Performing several integrations by parts and taking into account the boundary conditions $\eqref{pde-aug-system}_3-\eqref{pde-aug-system}_4$, the following holds: 
	\begin{gather} \label{proof-t1-2}
		\begin{aligned}[b]
			\int_{0}^{\ell} y^{(4)} \overline{v} \, dx
			=
			\Big[ y^{(3)} \overline{v} - y'' \overline{v'} \Big]_{x=0}^{x=\ell} + \int_{0}^{\ell} y'' \overline{v''} \, dx
			=
			\int_{0}^{\ell} y'' \overline{v''} \, dx;
		\end{aligned}
		\\
		\begin{aligned}
			\int_{0}^{1} z_s \overline{z}\, ds 
			=
			\big[|z|^2\big]_{s=0}^{s=1} - \int_{0}^{1} z \overline{z_s} \, ds 
			=
			|z(\cdot,1)|^2 - |v|^2  - \int_{0}^{1} z \nonumber \overline{z_s} \, ds
		\end{aligned}
		\\ \label{proof-t1-3}
		\Re \left( \int_{0}^{1} z_s \overline{z}\, ds \right)
		=
		\frac{1}{2} \Big( |z(\cdot,1)|^2 - |v|^2 \Big).
	\end{gather}
	Thus, combining \eqref{proof-t1-2} and \eqref{proof-t1-3}, we arrive at
	\begin{equation}
		\begin{aligned}[b]
			\langle \mathcal{A}U, U \rangle 
			&=
			- \mu_1 \int_{0}^{\ell} a(x) |v(x)|^2 \, dx 
			- \mu_2 \int_{0}^{\ell} a(x) z(x,1) \overline{v} \, dx  
			\\
			& \quad - \xi \tau^{-1} \int_{0}^{\ell} a(x) \int_{0}^{1} z_s \overline{z}\, ds \, dx + 2 i \Im \left( EI \int_{0}^{\ell} v'' \overline{y''} \, dx \right)
			\\
			\Re \langle \mathcal{A} U, U \rangle 
			&=
			- \mu_1 \int_{0}^{\ell} a(x) |v(x)|^2 \, dx  
			- \mu_2 \int_{0}^{\ell} a(x) \Re \left( z(x,1) \overline{v} \right) \, dx  
			\\
			& \quad - \xi \tau^{-1} \int_{0}^{\ell} a(x) \Re \left( \int_{0}^{1} z_s \overline{z}\, ds \right) \, dx
			\\
			&=
			\left( - \mu_1 + \frac{\xi \tau^{-1}}{2} \right) \int_{0}^{\ell} a(x) |v(x)|^2 \, dx 
			- \frac{\xi \tau^{-1}}{2} \int_{0}^{\ell} a(x) |z(x,1)|^2 \, dx 
			\\
			& \quad 
			- \mu_2 \int_{0}^{\ell} a(x) \Re \big( z(x,1) \overline{v(x)} \big) \, dx. \label{proof-t1-4}
		\end{aligned}
	\end{equation}
	Applying Young's inequality to the integrand of the last integral term, we get
	\begin{align}
		\Re \big( \mu_2 z(\cdot,1) \overline{v} \big)
		\leqslant
		\frac{|\mu_2|}{2} \left( |v|^2 + |z(\cdot,1)|^2 \right),
	\end{align}
	whose incorporation in \eqref{proof-t1-4} leads to:
	\begin{align*}
		\Re \langle \mathcal{A} U, U \rangle 
		&\leqslant
		\left( - \mu_1 + \frac{|\mu_2|}{2} + \frac{\xi \tau^{-1}}{2} \right) \int_{0}^{\ell} a(x) |v(x)|^2 \, dx
		\\
		& \quad 
		+ \left( \frac{|\mu_2|}{2} - \frac{\xi \tau^{-1}}{2} \right) \int_{0}^{\ell} a(x) |z(x,1)|^2 \, dx
		\\
		& \leqslant 
		- C \left( \int_{0}^{\ell} a(x) |v(x)|^2 \, dx + \int_{0}^{\ell} a(x) |z(x,1)|^2 \, dx \right),
	\end{align*}
	where the constant $C$ is positive by virtue of condition \eqref{aux-coef-cdt} and is defined as
	\[
	C := \min \left\{ 
	-\frac{|\mu_2|}{2} + \frac{\xi \tau^{-1}}{2}; \,
	\mu_1 - \frac{|\mu_2|}{2} - \frac{\xi \tau^{-1}}{2} 
	\right\}.
	\]
	Next, to prove the maximality of $\mathcal{A}$, it suffices to show that $\lambda I - \mathcal{A}$ is surjective for a given $\lambda > 0$. That is, for any $\Psi = (f,g,h)^\top \in \mathscr{H}$, we seek a solution $U = (y,v,z)^\top \in \mathscr{D}(\mathcal{A})$ to the equation $(\lambda I - \mathcal{A})U = \Psi$. Equivalently, this resolvent equation is expressed componentwise as: 
	\begin{gather} \label{proof-t1-5}
		\left\{
		\begin{aligned}
			& \lambda y - v = f \\
			& \lambda v + \frac{EI}{m} y^{(4)} + \frac{a}{m} \left( \mu_1 v + \mu_2 z(\cdot,1) \right) = g \\
			& \lambda z + \tau^{-1} z_s = h.
		\end{aligned} 
		\right.
	\end{gather}
	Solving equation $\eqref{proof-t1-5}_1$ and equation $\eqref{proof-t1-5}_3$ with respect to $x$ and $s$, respectively, and taking into account the boundary condition $\eqref{pde-aug-system}_4$, we obtain
	\begin{align}
		& v = \lambda y - f, \label{proof-t1-6}
		\\
		& z(\cdot,s) = e^{-\lambda \tau s} \left( v + \tau \int_{0}^{s} h(\cdot,r) e^{\tau \lambda r} \, dr \right). \label{proof-t1-7}
	\end{align}
	Consequently, inserting expressions \eqref{proof-t1-6}--\eqref{proof-t1-7} into $\eqref{proof-t1-5}_2$, we arrive at the following boundary value problem:
	\begin{gather} \label{proof-t1-8}
		\left\{
		\begin{aligned}
			& EI y^{(4)} + \lambda G_{\tau,1}(\lambda,\cdot) y = m g + G_{\tau,1}(\lambda,\cdot) f + G_{\tau}^h(\lambda,\cdot) \\
			& y(0) = y'(0) = y''(\ell) = y^{(3)}(\ell) = 0,
		\end{aligned} 
		\right.
	\end{gather}
	where the functions $G_{\tau,1}, \, G_\tau^h$ are defined by
	\begin{gather}
		G_{\tau,1}(\lambda,x) := m\lambda + \big( \mu_{1} + \mu_2 e^{-\tau \lambda} \big) a(x), 
		\\
		G_{\tau}^h(\lambda,x) := - m \mu_2 \tau e^{-\tau \lambda} \int_{0}^{1} h(\cdot,r) e^{\tau \lambda r} \, dr,
	\end{gather}
	We note that $\mu_1 + \mu_2 e^{-\tau \lambda} > 0$ by virtue of \eqref{aux-coef-cdt}. We multiply the equation $\eqref{proof-t1-8}_1$ by $\overline{\varphi}$ and integrate over $(0,\ell)$ with respect to $x$:
	\begin{align}
		\begin{split}
			& \int_0^\ell EI \left( y^{(4)} + \lambda G_{\tau,1}(\lambda,\cdot) y \right) \overline{\varphi} \, dx 
			\\
			& \quad 
			= \int_0^\ell  \left( m g + G_{\tau,1}(\lambda,\cdot) f + G_{\tau}^h(\lambda,\cdot)\right) \overline{\varphi} \, dx.
		\end{split}
	\end{align}
	By performing several integrations by parts together with the boundary conditions $\eqref{proof-t1-8}_2$, we obtain
	\begin{align}
		\begin{split}
			& EI \int_0^\ell y'' \overline{\varphi''} \, dx + \lambda \int_0^\ell G_{\tau,1}(\lambda,\cdot) \, y  \overline{\varphi} \, dx
			\\ 
			& \quad 
			= \int_0^\ell  \left( m g + G_{\tau,1}(\lambda,\cdot) f + G_{\tau}^h(\lambda,\cdot)\right) \overline{\varphi} \, dx.
		\end{split}
	\end{align}
	This leads to the variational problem which consists in finding $y \in H^2_e(0,\ell)$ such that
	\begin{align} \label{proof-t1-varpb}
		\mathcal{J}(y,\varphi) = \mathcal{L}(\varphi),
	\end{align}
	for all $\varphi \in H^2_e(0,\ell)$, where the functionals $\mathcal{J}$ and $\mathcal{L}$ are given by
	\begin{gather}
		\mathcal{J}(y,\varphi) 
		= EI \int_0^\ell y'' \overline{\varphi''} \, dx + \lambda  \int_0^\ell G_{\tau,1}(\lambda,\cdot) \, y \overline{\varphi} \, dx 
		\\
		\mathcal{L}(\varphi)
		= \int_0^\ell  \left( m g + G_{\tau,1}(\lambda,\cdot) f + G_{\tau}^h(\lambda,\cdot) \right) \overline{\varphi} \, dx.
	\end{gather}
	It is straightforward to verify that $\mathcal{J}$ and $\mathcal{L}$ are sesquilinear and antilinear, respectively.
	Furthermore, by applying the Cauchy--Schwarz inequality and taking into account condition \eqref{damper-cdt}, we obtain the following estimates:
	\begin{gather} \label{proof-t1-9}
		\left| \int_0^\ell y'' \overline{\varphi''} \, dx \right|
		\leqslant \lVert y'' \rVert_{L^2(0,\ell)} \lVert \varphi'' \rVert_{L^2(0,\ell)}
		\\ \label{proof-t1-10}
		\left| \int_0^\ell G_{\tau,1}(\lambda,\cdot) \, y \overline{\varphi} \, dx  \right|
		\leqslant 
		\Big[ m \lambda + \big( \mu_1 + \mu_2 e^{-\tau \lambda} \big) \lVert a \rVert_{L^\infty(0,\ell)} \Big] \lVert y \rVert_{L^2(0,\ell)} \lVert \varphi \rVert_{L^2(0,\ell)}.
	\end{gather}
	Combining expressions \eqref{proof-t1-9} and \eqref{proof-t1-10}, it follows that
	\begin{align}
		\big| \mathcal{J}(y,\varphi) \big| 
		& \leqslant 
		EI  \lVert y'' \rVert_{L^2(0,\ell)} \lVert \varphi'' \rVert_{L^2(0,\ell)} \nonumber 
		\\
		& \quad + \lambda \Big[ m \lambda + \big( \mu_1 + \mu_2 e^{-\tau \lambda} \big) \lVert a \rVert_{L^\infty(0,\ell)} \Big] \lVert y \rVert_{L^2(0,\ell)} \lVert \varphi \rVert_{L^2(0,\ell)} \nonumber
		\\
		\big| \mathcal{J}(y,\varphi) \big| 
		& \leqslant C \lVert y \rVert_{H^2(0,\ell)} \lVert \varphi \rVert_{H^2(0,\ell)}, \label{proof-t1-11}
	\end{align}
	where the positive constant $C$ is given by:
	\begin{equation}
		C := 2 \max \Big\{ 
		EI, \, \lambda \Big( m \lambda + \big( \mu_1 + \mu_2 e^{-\tau \lambda} \big) \lVert a \rVert_{L^\infty(0,\ell)} \Big).
	\end{equation} 
	Thus, $\mathcal{J}$ is continuous. Next, by the Cauchy--Schwarz inequality, we have
	\begin{align}
		\int_0^\ell \left( \int_{0}^{1} h(\cdot,r) e^{\tau \lambda r} \, dr \right) \,  \overline{\varphi} \, dx 
		& \leqslant
		\left\lVert \int_{0}^{1} h(\cdot,r) e^{\tau \lambda r} \, dr \right\rVert_{L^2(0,\ell)} \lVert \varphi \rVert_{L^2(0,\ell)} \nonumber
		\\ 
		& \leqslant
		e^{\tau \lambda} \lVert h \rVert_{L^2\big( (0,\ell) \times (0,1) \big)} \lVert \varphi \rVert_{L^2(0,\ell)}. \label{proof-t1-12}
	\end{align}
	Similarly to estimate \eqref{proof-t1-11}, using \eqref{proof-t1-12} yields
	\begin{align}
		\begin{split}
			\big| \mathcal{L}(\varphi) \big| 
			& \leqslant
			\Big[  
			m \lVert g \rVert_{L^2(0,\ell)}
			+ m \mu_2 \tau \lVert h \rVert_{L^2\big( (0,\ell) \times (0,1) \big)}
			\\
			& \quad 
			+ \Big( m \lambda + \big( \mu_1 + \mu_2 e^{-\tau \lambda} \big) \lVert a \rVert_{L^\infty(0,\ell)} \Big) \lVert f \rVert_{L^2(0,\ell)}  
			\Big] \lVert \varphi \rVert_{L^2(0,\ell)}
			\\
			& \leqslant 
			C \lVert \varphi \rVert_{H^2(0,\ell)},
		\end{split}
	\end{align}
	where the constant $C > 0$ is given by:
	\begin{equation}
		\begin{aligned}
			C 
			&:= 
			m \lVert g \rVert_{L^2(0,\ell)}
			+ m \mu_2 \tau \lVert h \rVert_{L^2( (0,\ell) \times (0,1) )} 
			\\
			& \quad 
			+ \Big( m \lambda + \big( \mu_1 + \mu_2 e^{-\tau \lambda} \big) \lVert a \rVert_{L^\infty(0,\ell)} \Big) \lVert f \rVert_{L^2(0,\ell)}. 
		\end{aligned}
	\end{equation} 
	This proves the continuity of the functional $\mathcal{L}$. Moreover, $\mathcal{J}$ is coercive on $H^2_e(0,\ell)$ since
	\begin{align}
		\Re \mathcal{J}(y,y) 
		& = EI \int_0^\ell \big| y'' \big|^2 \, dx + \lambda \int_0^\ell G_{\tau,1}(\lambda,x) \, |y|^2 \, dx \nonumber
		\\
		& \geqslant EI \int_0^\ell \big| y'' \big|^2 \, dx + m \lambda^2 \int_0^\ell |y|^2 \, dx + (1- e^{-\tau \lambda}) \mu_1 \lambda \int_0^\ell a(x) |y|^2 \, dx \nonumber
		\\
		\Re \mathcal{J}(y,y) 
		& \geqslant \min \left\{ EI; m \lambda^2 \right\} |y|^2_2 \geqslant \min \left\{ EI; m \lambda^2 \right\} C_1 \lVert y \rVert^2_{H^2(0,\ell)}.
	\end{align}
	Consequently, by the Lax--Milgram theorem, the variational problem \eqref{proof-t1-varpb} admits a unique solution $y \in H^2_e(0,\ell)$.
	
	Conversely, consider a weak solution $y \in H^2_e(0,\ell)$. Applying integration by parts, equation \eqref{proof-t1-varpb} becomes:
	\begin{align}
		& \int_0^\ell \left( EI y^{(4)} + \lambda G_{\tau,1}(\lambda,\cdot) \, y  \right) \overline{\varphi} \, dx 
		+ EI y''(\ell) \overline{\varphi'(\ell)} - EI y^{(3)} (\ell) \overline{\varphi(\ell)} \nonumber
		\\
		& \quad 
		= \int_0^\ell  \left( m g + G_{\tau,1}(\lambda,\cdot) f + G_{\tau}^h(\lambda,\cdot)\right) \overline{\varphi} \, dx. \label{proof-t1-13}
	\end{align}
	Taking first $\varphi \in \mathscr{D}(0,\ell)$, the preceding identity reduces to
	\begin{equation}
		\begin{gathered}
			\int_0^\ell \left( EI y^{(4)} + \lambda G_{\tau,1}(\lambda,\cdot) \, y  \right) \overline{\varphi} \, dx 
			\int_0^\ell  \left( m g + G_{\tau,1}(\lambda,\cdot) f + G_{\tau}^h(\lambda,\cdot)\right) \overline{\varphi} \, dx.
		\end{gathered}
	\end{equation}
	This implies that $EI y^{(4)} + \lambda G_{\tau,1}(\lambda,\cdot) \, y = m g + G_{\tau,1}(\lambda,\cdot) f + G_{\tau,2}(\lambda)$ almost everywhere on $(0,\ell)$. Consequently, $y \in H^4(0,\ell)$ and satisfies $\eqref{proof-t1-8}_1$ by virtue of the density of $\mathscr{D}(0,\ell)$ in $L^2(0,\ell)$.
	
	Returning to $\varphi \in H^2_E(0,\ell)$ and from the above, equation \eqref{proof-t1-13} leads to:
	\begin{equation}
		EI y''(\ell) \overline{\varphi'(\ell)} - EI y^{(3)} (\ell) \overline{\varphi(\ell)} = 0.
	\end{equation}
	By the arbitrariness of $\varphi$, $y$ fulfills the boundary condition $\eqref{proof-t1-8}_2$, and therefore solves \eqref{proof-t1-8}.
	Thus, the operator $\lambda I - \mathcal{A}$ is surjective, which proves that $\mathcal{A}$ is maximal.
	
	Finally, the Lumer--Phillips theorem ensures that $\mathcal{A}$ generates a $C_0$-semigroup of contractions on $\mathscr{H}$, which guarantees the existence and uniqueness of solutions to the Cauchy problem \eqref{Cauchy-pb}.
\end{proof}

\section{Exponential stability} \label{sec3}

This section is dedicated to the statement and proof of the exponential stability of system \eqref{pde-system} based on Lyapunov method combined with multiplier techniques. To this end, we first establish several technical lemmas.

We begin by introducing the auxiliary functional $Q(t)$, defined as
\begin{gather} \label{def:aux-funct}
	Q(t) := Q_1(t) - Q_2(t), \quad t \geqslant 0,
\end{gather}
where the functionals $Q_i, \, i=1,2$ are given by
\begin{gather}
	Q_1(t) := m \int_{0}^{\ell} y y_t \, dx
	\\
	Q_2(t) := \int_{0}^{\ell} a(x) \int_{0}^{1} s y_t^2(x,t-\tau s) \, ds \, dx.
\end{gather}
Prior to estimating the time derivative of $Q(t)$, we state and prove a few Poincaré-type inequalities.
\begin{lemma} \label{l2}
	Let $y \in H^2_e(0,\ell)$. Then, the following estimates hold:
	\begin{gather}
		\int_0^\ell y^2(x) \, dx 
		\leqslant
		\frac{\ell^4}{8} \big\lVert y'' \big\rVert_{L^2(0,\ell)}^2, \label{ineq-1:Poincare}
		\\
		\int_0^\ell \, a(x) \, y^2(x) \, dx 
		\leqslant 
		\frac{\ell^3}{2} \int_0^\ell \, a(x) \, dx \, \big\lVert y'' \big\rVert_{L^2(0,\ell)}^2 \label{ineq-2:Poincare}.
	\end{gather}
\end{lemma}

\begin{proof}
	We have for all $x \in (0,\ell)$:
	\begin{gather}
		y^2(x) 
		= \left( \int_0^x \, y'(r) \, dr \right)^2 
		\leqslant 
		x \int_0^x \, y'(r)^2 \, dr. \label{proof-l2-1}
	\end{gather}
	In the same manner, taking into account $y'(0) = 0$, we have
	\begin{gather}
		y'^2(r) 
		\leqslant 
		r \int_0^r \, y''(\sigma)^2 \, d\sigma 
		\leqslant 
		r \big\lVert y'' \big\rVert_{L^2(0,\ell)}^2. \label{proof-l2-2}
	\end{gather}
	Consequently, integrating over $(0,\ell)$ and combining estimates \eqref{proof-l2-1} and \eqref{proof-l2-2}, we deduce that
	\begin{gather*}
		\begin{aligned}[b]
			\int_0^\ell y^2(x) \, dx 
			& \leqslant 
			\int_0^\ell x \int_0^x \, y'(r)^2 \, dr \, dx
			\\
			& \leqslant 
			\int_0^\ell x \int_0^x \, r \, dr \, dx \; \big\lVert y'' \big\rVert_{L^2(0,\ell)}^2 
			\\
			\int_0^\ell y^2(x) \, dx
			& \leqslant
			\frac{\ell^4}{8} \, \big\lVert y'' \big\rVert_{L^2(0,\ell)}^2;
		\end{aligned}
		\\
		\begin{aligned}[b]
			\int_0^\ell a(x) y^2(x) \, dx 
			& \leqslant 
			\int_0^\ell x a(x) \int_0^x \, r \, dr \, dx \, \big\lVert y'' \big\rVert_{L^2(0,\ell)}^2 
			\\
			& \leqslant
			\int_0^\ell \frac{x^3}{2} a(x) \, dx \, \big\lVert y'' \big\rVert_{L^2(0,\ell)}^2
			\\
			\int_0^\ell a(x) y^2(x) \, dx
			& \leqslant
			\frac{\ell^3}{2} \int_0^\ell  a(x) \, dx \, \big\lVert y'' \big\rVert_{L^2(0,\ell)}^2.
		\end{aligned}
	\end{gather*}
	This completes the proof.
\end{proof}

\begin{lemma} \label{l3}
	Let $y$ be the solution of system \eqref{pde-system}. For every $\varepsilon_1>0$, the time derivative of $Q(t)$ satisfies for all $t \geqslant 0$,
	\begin{gather} \label{estim:aux-fct-derivative}
		\begin{aligned}[b]
			Q'(t)
			& \leqslant 
			\int_{0}^{\ell} \left( \frac{m}{2} + \frac{\mu_1^2}{2 \varepsilon_1} a(x) \right) y_t^2 \, dx 
			+ \int_{0}^{\ell} \left( -EI + \varepsilon_1 \frac{\ell^3}{2}  \int_{0}^{\ell} a(x) \, dx \right) y_{xx}^2 \, dx 
			\\ 
			& \quad 
			+ \left( \tau^{-1} + \frac{\mu_2^2}{2 \varepsilon_1} \right) \int_{0}^{\ell} a(x) y_t^2(x,t-\tau) \, dx
			- \tau^{-1} \int_{0}^{\ell} a(x) \int_{0}^{1} y_{t}^2(x,t-\tau s) \, ds \, dx.
		\end{aligned}
	\end{gather}
\end{lemma}

\begin{proof}
	First, we compute the time derivative of $Q_1(t)$. Using equation $\eqref{pde-system}_1$, we have
	\begin{gather} \label{proof-l3-1}
		\begin{aligned}[b]
			Q_1'(t) 
			& = 
			m \int_{0}^{\ell} \big( y_t^2 + y y_{tt} \big) \, dx 
			\\
			& =
			m \int_{0}^{\ell} y_t^2 \, dx 
			- EI \int_{0}^{\ell} y y_{xxxx} \, dx 
			- \int_{0}^{\ell} a(x) \Big( \mu_1 y_t + \mu_2 y_t(x,t-\tau) \Big) y \, dx.
		\end{aligned}
	\end{gather}
	Next, using integrations by parts together with the boundary condition $\eqref{pde-system}_2$, the second integral term in \eqref{proof-l3-1} can be rewritten as
	\begin{equation} 
		\begin{aligned}[b]
			EI \int_{0}^{\ell} y_{xxxx} y \, dx
			& = 
			EI \big[ y y_{xxx} - y_x y_{xx} \big]_{x=0}^{x=\ell} + EI \int_{0}^{\ell} y_{xx}^2 \, dx
			\\
			& = EI \int_{0}^{\ell} y_{xx}^2 \, dx.
		\end{aligned}
	\end{equation}
	It follows that: 
	\begin{equation}  \label{proof-l3-2}
		\begin{split}
			Q_1'(t) 
			& = 
			m \int_{0}^{\ell} y_t^2 \, dx 
			- EI \int_{0}^{\ell} y_{xx}^2 \, dx 
			- \int_{0}^{\ell} a(x) \Big( \mu_1 y_t + \mu_2 y_t(x,t-\tau) \Big) y \, dx.
		\end{split}
	\end{equation}
	Furthermore, applying Young's inequality to the last term of \eqref{proof-l3-2} and invoking estimate \eqref{ineq-2:Poincare} from Lemma \ref{l2}, we obtain for any $\varepsilon_1 > 0$:
	\begin{gather} \label{proof-l3-3}
		\begin{aligned}[b]
			& \left| \int_{0}^{\ell} x a(x) \Big( \mu_1 y_t + \mu_2 y_t(x,t-\tau) \Big) y \, dx \right| 
			\\ 
			& \quad \leqslant 
			\varepsilon_1 \int_{0}^{\ell} a(x) y^2 \, dx  
			+ \frac{\mu_1^2}{2 \varepsilon_1} \int_{0}^{\ell} a(x) y_t^2 \, dx 
			+ \frac{\mu_2^2}{2 \varepsilon_1} \int_{0}^{\ell} a(x) y_t^2(x,t-\tau) \, dx 
			\\
			& \quad \leqslant 
			\frac{\ell^3}{2} \varepsilon_1 \int_0^\ell \, a(x) \, dx \, \big\lVert y'' \big\rVert_{L^2(0,\ell)}^2
			+ \frac{\mu_1^2}{2 \varepsilon_1} \int_{0}^{\ell} a(x) y_t^2 \, dx 
			+ \frac{\mu_2^2}{2 \varepsilon_1} \int_{0}^{\ell} a(x) y_t^2(x,t-\tau) \, dx.
		\end{aligned}
	\end{gather}
	Thus, inserting estimate \eqref{proof-l3-3} into \eqref{proof-l3-2} leads to
	\begin{gather} \label{proof-l3-4}
		\begin{aligned}[b]
			Q_1'(t) 
			& \leqslant 
			\int_{0}^{\ell} \left( m + \frac{\mu_1^2}{2 \varepsilon_1} a(x) \right)y_t^2 \, dx 
			\\
			& \quad 
			+ \left( - EI + \varepsilon_1 \frac{\ell^3}{2} \int_0^\ell \, a(x) \, dx \right) \int_{0}^{\ell} y_{xx}^2 \, dx 
			\\ 
			& \quad 
			+ \frac{\mu_2^2}{2 \varepsilon_1} \int_{0}^{\ell} a(x) y_t^2(x,t-\tau) \, dx.
		\end{aligned}
	\end{gather}
	Next, differentiating $Q_2(t)$ with respect to $t$, we obtain
	\begin{gather} \label{proof-l3-5}
		\begin{aligned}[b]
			Q_2'(t) 
			& = 2 \int_{0}^{\ell} a(x) \int_{0}^{1} s y_t(x,t-\tau s) y_{tt}(x,t-\tau s) \, ds \, dx 
			\\
			& = 
			- \tau^{-1} \int_{0}^{\ell} a(x) \int_{0}^{1} s \mathcolor{blue}{\frac{\partial}{\partial s}} \big( y_{t}^2(x,t-\tau s) \big) \, ds \, dx 
			\\
			& = 
			\tau^{-1} \int_{0}^{\ell} a(x) \int_{0}^{1} y_{t}^2(x,t-\tau s) \, ds \, dx 
			- \tau^{-1}  \int_{0}^{\ell} a(x) y_{t}^2(x,t-\tau) \, dx.
		\end{aligned}
	\end{gather}
	Finally, combining estimates \eqref{proof-l3-4} and \eqref{proof-l3-5} directly yields \eqref{estim:aux-fct-derivative}.
\end{proof}

We now define the Lyapunov functional $V(t)$ by 
\begin{align}
	V(t) := E(t) + \delta Q(t), \quad t \geqslant 0,
\end{align} 
where $E(t)$, $Q(t)$ are defined in \eqref{energy-system} and \eqref{def:aux-funct}, respectively, and $\delta>0$ is a positive parameter to be specified later.
\begin{lemma} \label{l4}
	For a sufficiently small $\delta>0$, there exist two positive constants $\theta_1$ and $\theta_2$ such that, for all $t \geqslant 0$:
	\begin{gather} \label{equiv-energy-lyap}
		\theta_1 E(t) \leqslant V(t) \leqslant \theta_2 E(t).
	\end{gather}
\end{lemma}

\begin{proof}
	Using estimate \eqref{ineq-1:Poincare} together with Young's inequality, we have for any $\varepsilon_2 > 0$:
	\begin{equation}
		\begin{aligned}[b]
			|Q_1(t)| 
			& \leqslant
			m \int_{0}^{\ell} |y y_t| \, dx
			\\
			& \leqslant
			m \left( \frac{\varepsilon_2}{2} \int_{0}^{\ell} y^2 \, dx + \frac{1}{2\varepsilon_2} \int_{0}^{\ell} y_t^2 \, dx \right)
			\\
			& \leqslant 
			m \left( \frac{\varepsilon_2 \ell^4}{16}  \lVert y_{xx} \rVert_{L^2(0,\ell)}^2 + \frac{1}{2\varepsilon_2} \lVert y_t \rVert_{L^2(0,\ell)}^2 \right).
		\end{aligned}
	\end{equation}
	Choosing $\varepsilon_2 = \frac{2\sqrt{2}}{\ell^2}$, the preceding estimate simplifies to:
	\begin{align}
		|Q_1(t)| \leqslant \frac{m \ell^2}{4 \sqrt{2}} \left( \lVert y_{xx} \rVert_{L^2(0,\ell)}^2 + \lVert y_t \rVert_{L^2(0,\ell)}^2 \right). \label{proof-l4-1}
	\end{align}
	On the other hand, $Q_2(t)$ satisfies:
	\begin{equation}
		\begin{aligned}[b]
			|Q_2(t)| 
			& \leqslant
			\int_{0}^{\ell} a(x) \int_{0}^{1} y_t^2(x,t-\tau s) \, ds \, dx. \label{proof-l4-2}
		\end{aligned}
	\end{equation}
	Gathering \eqref{proof-l4-1} and \eqref{proof-l4-2} leads to
	\begin{equation}
		\begin{aligned}[b]
			|Q(t)| 
			& \leqslant 
			\frac{m \ell^2}{4 \sqrt{2}} \left( \lVert y_{xx} \rVert_{L^2(0,\ell)}^2 + \lVert y_t \rVert_{L^2(0,\ell)}^2 \right)
			\\
			& \quad 
			+ \delta \int_{0}^{\ell} a(x) \int_{0}^{1} y_t^2(x,t-\tau s) \, ds \, dx
			\\
			& \leqslant C E(t), \label{proof-l4-3}
		\end{aligned}
	\end{equation}
	where $C > 0$ is defined as
	\begin{align}
		C := \max \left\{ \frac{2 \delta}{\xi}, \, \frac{\ell^2}{2 \sqrt{2}}, \, \frac{m \ell^2}{2 \sqrt{2} EI} \right\}.
	\end{align}
	Finally, by applying the triangle inequality on the estimate \eqref{proof-l4-3} and by adding the energy $E(t)$ to each side of the resulting inequality, we arrive at: 
	\begin{gather*}
		0 < (1 - \delta C) E(t) \leqslant V(t) \leqslant (1 + \delta C) E(t),
	\end{gather*} 
	by choosing $\delta \in (0, C^{-1})$.
\end{proof}

The main result of this paper states as follows:
\begin{theorem} \label{t2}
	Suppose that assumptions \eqref{internal-coef-cdt} and \eqref{aux-coef-cdt} are satisfied. For any regular solution $y$ of \eqref{pde-system}, there exist constants $\upsilon_1, \upsilon_2 > 0$ such that the energy $E(t)$ satisfies
	\begin{align} \label{estim:expo-decay}
		E(t) \leqslant \upsilon_1 E(0) e^{-\upsilon_2 t}, \quad \forall t \geqslant 0.
	\end{align}
\end{theorem}

\begin{proof}
	In view of Lemmas \ref{l1a} and \ref{l3}, we obtain \eqref{decrease-energy} and \eqref{estim:aux-fct-derivative}, together with condition $\eqref{damper-cdt}_2$, the time derivative of $V(t)$ satisfies
	\begin{gather} \label{proof-t2-1}
		\begin{aligned}[b]
			& \frac{d}{dt} V(t) 
			= \frac{d}{dt} E(t) + \frac{d}{dt} Q(t) 
			\\
			& \leqslant 
			\int_{0}^{\ell} \bigg[ -\left( -\mu_1 + \frac{|\mu_2|}{2} + \frac{\xi \tau^{-1}}{2} \right) a(x) 
			+ \delta \left( \frac{m}{2} + \frac{\mu_1^2}{2 \varepsilon_1} a(x) \right) \bigg] y_t^2 \, dx 
			\\
			& \quad 
			+\delta \int_{0}^{\ell} \left( -EI + \varepsilon_1 \frac{\ell^3}{2}  \int_{0}^{\ell} a(x) \, dx \right) y_{xx}^2 \, dx 
			\\ 
			& \quad 
			- \delta \tau^{-1} \int_{0}^{\ell} a(x) \int_{0}^{1} y_{t}^2(x,t-\tau s) \, ds \, dx 
			\\
			& \quad 
			+ \left[ \frac{|\mu_2|}{2} - \frac{\xi \tau^{-1} }{2} + \delta \left( \tau^{-1} + \frac{\mu_2^2}{2 \varepsilon_1} \right) \right] \int_{0}^{\ell} a(x) y_t^2(x,t-\tau) \, dx
			\\
			& \leqslant 
			\bigg[ -\left( -\mu_1 + \frac{|\mu_2|}{2} + \frac{\xi \tau^{-1}}{2} \right) a_0 
			+ \delta \left( \frac{m}{2} + \frac{\mu_1^2}{2 \varepsilon_1} \lVert a \rVert_{L^\infty(0,\ell)} \right) \bigg] \int_{0}^{\ell} y_t^2 \, dx 
			\\
			& \quad 
			+\delta \int_{0}^{\ell} \left( -EI + \varepsilon_1 \frac{\ell^3}{2}  \int_{0}^{\ell} a(x) \, dx \right) y_{xx}^2 \, dx 
			\\ 
			& \quad 
			- \delta \tau^{-1} \int_{0}^{\ell} a(x) \int_{0}^{1} y_{t}^2(x,t-\tau s) \, ds \, dx 
			\\
			& \quad 
			+ \left[ \frac{|\mu_2|}{2} - \frac{\xi \tau^{-1} }{2} + \delta \left( \tau^{-1} + \frac{\mu_2^2}{2 \varepsilon_1} \right) \right] \int_{0}^{\ell} a(x) y_t^2(x,t-\tau) \, dx.
		\end{aligned}
	\end{gather} 
	We fix $\varepsilon_1 > 0$ as
	\begin{equation} \label{epsilon1}
		\varepsilon_1 = EI \left( \ell^3 \int_{0}^{\ell} a(x) \, dx \right)^{-1}.
	\end{equation}
	Then, we select $\delta$ small enough so that $0 < \delta < \delta_0$, where $\delta_0$ is defined by
	\begin{gather}
		\begin{aligned}[b]
			\delta_0 
			& := \min \Bigg\{ 
			\left( - \frac{|\mu_2|}{2} + \frac{\xi}{2 \tau} \right) \left( \tau^{-1} + \frac{\mu_2^2}{2 \varepsilon_1} \right)^{-1},
			\\
			& \qquad \qquad 
			\left( -\mu_1 + \frac{|\mu_2|}{2} + \frac{\xi}{2 \tau} \right) \left( \frac{m}{2} + \frac{\mu_1^2}{2 \varepsilon_1} \lVert a \rVert_{L^\infty(0,\ell)} \right)^{-1} 
			\Bigg\}.
		\end{aligned}
	\end{gather}
	we deduce that
	\begin{gather}
		V'(t) \leqslant - C_\delta E(t),
	\end{gather} 
	where the positive constant $C_\delta$ is given by:
	\begin{gather*}
		C_\delta:= \min \left\{ \delta, \, \frac{2 \delta}{\tau \xi}, \, 
		\left( -2\mu_1 + |\mu_2| + \frac{\xi}{2 \tau} \right) \frac{a_0}{m} - \delta \left( 1 + \frac{\mu_1^2}{\varepsilon_1 m} \lVert a \rVert_{L^\infty(0,\ell)} \right) \right\}.
	\end{gather*}
	Thanks to the functional equivalence established in Lemma \ref{l4}, we arrive at
	\begin{gather}
		V'(t) \leqslant - \frac{C_\delta}{\theta_1} V(t).
	\end{gather}
	applying Gronwall's inequality together with Lemma \ref{l4} directly yields the exponential decay estimate \eqref{estim:expo-decay}.
\end{proof}
	

\end{document}